\documentclass[12pt]{amsart}
\usepackage{amssymb,amscd,amsthm,amsmath,amsfonts,epsf,latexsym}
\usepackage{graphicx}
\usepackage{epsfig}

\newtheorem{thm}{Theorem}[section]
\newtheorem{theorem}[thm]{Theorem}

\newtheorem{lemma}[thm]{Lemma}

\newtheorem{prop}[thm]{Proposition}

\theoremstyle{remark}
\newtheorem{remark}[thm]{Remark}

\theoremstyle{definition}
\newtheorem{defn}[thm]{Definition}

\numberwithin{equation}{section}

\newcommand{\Hk}{\mathcal{H}}

\newcommand{\C}{{\mathbb{C}}}

\newcommand{\Z}{{\mathbb{Z}}}

\renewcommand{\O}{{\mathcal O}}

\newcommand{\tr}{\mathrm{tr}}
\newcommand{\sgn}{\mathrm{sgn}}
\newcommand{\Tr}{\mathrm{Tr}}
\newcommand{\GL}{\mathrm{GL}}

\newcommand{\End}{\mathrm{End}}

\newcommand{\PSU}{\mathrm{PSU}}

\newcommand{\ST}{\s^2T_m}

\newcommand{\PSp}{\mathrm{PSp}}
\newcommand{\PGL}{\mathrm{PGL}}

\newcommand{\TL}{\overline{T}_{m}}

\newcommand{\I}{\mathcal{I}}
\newcommand{\s}{\mathcal{S}}

\newcommand{\B}{\mathcal{B}}
\newcommand{\CC}{\mathcal{C}}

\newcommand{\oV}{\overline{V}}
\newcommand{\oW}{\overline{W}}

\newcommand{\T}{\mathcal{T}}
\newcommand{\tla}{\tilde{\la}}

\newcommand{\Ga}{\Gamma}

\newcommand{\N}{\mathbb{N}}

\newcommand{\lra}{\leftrightarrow}

\newcommand{\tG}{\tilde{G}}
\newcommand{\tE}{\tilde{E}}
\newcommand{\la}{{\lambda}}
\newcommand{\La}{{\Lambda}}

\begin{document}
\title[An algebra-level link-polynomial identity]
{An algebra-level version of a link-polynomial identity of Lickorish}

\author{Michael J. Larsen}
\email{larsen@math.indiana.edu}
\address{Department of Mathematics\\
    Indiana University \\
    Bloomington, IN 47405\\
    U.S.A.}

\author{Eric C. Rowell}
\email{errowell@indiana.edu}
\address{Department of Mathematics\\
    Indiana University \\
    Bloomington, IN 47405\\
    U.S.A.}

\thanks{The authors are partially supported by NSF grant  DMS-034772.}

\begin{abstract}
We establish isomorphisms between certain specializations of BMW
algebras and the symmetric squares of Temperley-Lieb algebras.
These isomorphisms imply a link-polynomial
identity due to W.~B.~R.~Lickorish. As an application, we compute the closed images
of the irreducible braid group representations factoring over these
specialized BMW algebras.
\end{abstract}
\maketitle
\section{Introduction}

In \cite{Li}, W.~B.~R.~Lickorish proved the following relation between
values of the Kauffman and Jones polynomials of an oriented link:
\begin{equation}
F_L(q^3,q^{-1}+q) = (-1)^{c(L)-1}V_L(-q^{-2}).
\label{Lickorish}
\end{equation}
This identity turns out to be a manifestation of a broader
phenomenon.  There exist two families of finite dimensional algebras
(actually von Neumann algebras): on the one hand, Birman-Murakami-Wenzl algebras with
a relation between the two parameters suggested by (\ref{Lickorish}),
and on the other, symmetric squares of Temperley-Lieb algebras.
On each side we have a natural trace
and a natural homomorphism from the group algebra of a braid group.
We show that there is a natural
isomorphism between corresponding algebras which respects both
structures and therefore ``explains'' (\ref{Lickorish}).  The
equality of dimensions gives a new combinatorial identity which can be expressed as an
explicit bijection between ``oscillating'' Young tableaux and pairs of ordinary tableaux.
Interestingly, our proof of the algebra isomorphism depends on first
establishing the combinatorial result; this allows us to show that the natural homomorphism
is actually an isomorphism.

The original motivation for this paper was our attempt to understand
the closed images of braid groups in the (projective) unitary
representations associated with the Kauffman polynomial at $q=e^{\pi
i/\ell}$. (For the HOMFLY polynomial, this was done by the
first-named author together with M.~Freedman and Z.~Wang
\cite{FLW}.)   A preliminary analysis of the Kauffman polynomial
case was undertaken by both authors and Wang \cite{LRW}.   The image
of any half-twist has eigenvalues $q$, $-q^{-1}$, and $r^{-1}$, and
\cite{LRW} excludes the cases where the ratio of two eigenvalues is
$-1$, or where all three eigenvalues lie in geometric progression.
Certain cases (for example $r=q^{-1}$) degenerate to those
considered in \cite{FLW}, and the cases excluded in \cite{LRW} that
remain are
 $r=q$, $r=\pm i$, and $r=q^{3}$. The first set of exceptions
will be discussed in the doctoral dissertation of Jennifer Franko;
in these cases, the image groups are finite.  The second set of
exceptions appears to be connected with self-dualities.  This paper
arose from our discovery that the third set of exceptions had a clean
algebraic interpretation. The actual classification of closed images
is given in Theorem~\ref{images} in the final section of the paper.

\subsection*{Acknowledgements} E. Rowell would like to thank Vaughan
Jones for useful conversations.
\section{Combinatorial Notation and Results}\label{comb}
The combinatorial language of Young diagrams plays a major role in
what follows so we establish notation and terminology for later use.

A Young diagram $\la$ is an array of boxes so that the number of
boxes in each row (resp. column) decreases weakly as one reads
downwards (resp. to the right), and we denote the set of Young
diagrams by $YD$.  Denote by $\la_i$ (resp. $\tla_i$) the number of
boxes in the $i$th row (resp. column) of $\la$.
 We identify $\la$ with an ordered list of its rows
$\la=[\la_1,\la_2,\ldots,\la_k]$ or columns
$\la:=[\tla_1,\ldots,\tla_j]^t$. The size $|\la|$ is defined to be
the total number of boxes
$|\la|=\la_1+\la_2+\cdots+\la_k=\tla_1+\cdots\tla_j$. If
$\la_i\leq\mu_i$ for all $i$ (where $\la_i=0$ is permitted) we write
$\la\subset\mu$, and if in addition $\mu$ can be obtained from $\la$
by adding one box, we write $\la\rightarrow\mu$.  The relation
$\subset$ is encoded in \emph{Young's lattice}.  An increasing path
in Young's lattice from $[0]$ to $\la$
$$t_\la:[0]=\la^{(0)}\rightarrow\la^{(1)}\rightarrow\cdots\rightarrow\la^{(m)}=\la$$
 where $|\la^{(j)}|=j$ is called a \emph{Young tableau of shape
$\la$}.  Denote by $\T(\la)$ the set of Young tableaux of shape
$\la$. We shall be particularly interested in Young tableaux $t_\la$
whose shapes $\la^{(j)}$ are restricted to a subset of $YD$.  In
particular we define a set $$\La(j,\ell):=\{[j-p,p],0\leq
j-2p\leq\ell-2\}$$ consisting of Young diagrams of size $j$ with at
most 2 rows whose row-difference is bounded by $\ell-2$. This
definition makes sense for $3\leq\ell\leq\infty$, where the case
$\ell=\infty$ corresponds to the set of \emph{all} Young diagrams of
size $j$ with at most 2 rows.  Notice that $\la\in\La(j,\ell)$ is
completely determined by its first row, $\la_1$, since
$\la_2=j-\la_1$.  We set $\La(\ell)=\bigcup_{0\leq j}\La(j,\ell)$ so
that $\La(\infty)$ is the subset of $YD$ consisting of all diagrams
with at most 2 rows.  Then we denote by $\T_\ell(\la)$ the set of
all restricted Young tableaux $t_\la$ where each
$\la^{(j)}\in\La(j,\ell)$.  Observe that $\T_\infty(\la)=\T(\la)$
since any Young tableaux terminating at a diagram
$\la\in\La(\infty)$ can only pass through diagrams in $\La(\infty)$.

These notions can be generalized: if $\la\rightarrow\mu$ \emph{or}
$\mu\rightarrow\la$, \emph{i.e.} $\la$ and $\mu$ differ by one box,
we write $\la\lra\mu$. A general path of length $m$ from $[0]$ to
$\la$ in Young's lattice
$$o_\la:[0]=\la^{(0)}\lra\la^{(1)}\lra\cdots\lra\la^{(m)}=\la$$ is
called an \emph{oscillating tableau of length $m$ and shape $\la$}.
Observe that $j-|\la^{(j)}|$ is always a non-negative even number.
We denote by $\O(m,\la)$ the set of oscillating tableaux of length
$m$ and shape $\la$.  We will often restrict the shapes to a subset
of $YD$, in this case the set:
$$\Ga(\ell):=\{\la\in YD:\tla_1+\tla_2\leq 4, \la_1+\la_2\leq\ell-2\}\cup\{[\ell-2,1^2]\},$$
where we will be interested in the (non-degenerate) cases:
$6\leq\ell\leq\infty$. For $\ell=\infty$ the conditions reduce to
$\tla_1+\tla_2\leq 4$. We denote by $\O_\ell(m,\la)$ the set of
oscillating tableaux of length $m$ and shape $\la$ restricted to the
set $\Ga(\ell)$.

Our basic combinatorial result is:
\begin{theorem}\label{maincomb}
Let $\la,\mu\in\La(m,\ell)$, and define $\nu_1=\la_1+\mu_1-m$ and
$\nu_2=|\la_1-\mu_1|$.  Then if $\la\neq\mu$ we have:
\begin{equation}
\label{tensordim} |\O_\ell(m,[\nu_1,\nu_2])| =
|\T_\ell(\la)|\cdot|\T_\ell(\mu)| ,
\end{equation}
while if $\la=\mu$ we have:
\begin{equation}
\label{symdim} |\O_\ell(m,[\nu_1])| =  \binom{|\T_\ell(\la)|+1}{2}
\end{equation}
and
\begin{equation}
\label{altdim} |\O_\ell(m,[\nu_1]^*)| = \binom{|\T_\ell(\la)|}{2},
\end{equation}
where $\ast$ is defined in (\ref{stardef}) below.
\end{theorem}
\begin{proof}

For each $m\ge 1$ and $\ell\ge 6$, we construct an explicit
bijection between two kinds of objects: on one side, pairs of
restricted tableaux $(t_\lambda,t_\mu)$, where
$\lambda,\mu\in\Lambda(m,\ell)$ and $\lambda_1\ge\mu_1$ and on the
other, oscillating tableaux $o_\nu$ of length $m$ with shapes
restricted to $\Gamma(\ell)$.

On $\Ga(\ell)$ define a reflection $\ast$ by:
\begin{equation}\label{stardef}
\text{$\tla^*_1=4-\tla_1$ and $\tla^*_j=\tla_j$ for
$j>1$}\end{equation} Explicitly, we have $[0]^*=[1^4]$,
$[\la_1]^*=[\la_1,1^2]$ for $\la_1>0$ and
$[\la_1,\la_2]^*=[\la_1,\la_2]$.

For any $(\sigma,\tau)\in\La(m,\ell)\times\La(m,\ell)$ define the
following functions:
 \begin{enumerate}
 \item $f(\sigma,\tau):=\sigma_1+\tau_1-m$
 \item $g(\sigma,\tau):=|\sigma_1-\tau_1|$
 \item $s(\sigma,\tau):=\sgn(\sigma_1-\tau_1)$
 \end{enumerate}

Suppose we are given $\la,\mu\in\La(m,\ell)$ with $\la_1\geq\mu_1$,
and $(t_\la,t_\mu)\in\T_\ell(\la)\times\T_\ell(\mu)$, \emph{i.e.}
$$t_\la:=\la^{(0)}=[0]\rightarrow\la^{(1)}\rightarrow\cdots\rightarrow\la^{(m)}=\la$$
and
$$t_\mu:=\mu^{(0)}=[0]\rightarrow\mu^{(1)}\rightarrow\cdots\rightarrow\mu^{(m)}=\mu.$$
The basic idea is that the two rows of the $j$th term in the
oscillating tableau associated to $\lambda$ and $\mu$ are obtained
by plugging $\lambda^{(j)}$ and $\mu^{(j)}$ into the formulas (1)
and (2) above.  Sometimes, however, the resulting diagram must be
reflected.  We now explain the rules for determining when this must
be done.

For each $j$, let $m_j$ denote the maximal positive integer $i\le j$
such that $s(\lambda^{(i)},\mu^{(i)})\neq 0$; if no such integer
exists, let $m_j=0$.  We define
\begin{equation*}
s^{(j)} =
\begin{cases}
1&\text{if $m_j=0$,}\\
s(\lambda^{(m_j)},\mu^{(m_j)})& \text{if $m_j\neq 0$.}
\end{cases}
\end{equation*}
We construct an oscillating tableau $o_\nu$ as follows.  For each
$j$ define
\begin{equation*}
\nu^{(j)}=
\begin{cases}
[f(\la^{(j)},\mu^{(j)}),g(\la^{(j)},\mu^{(j)})],&\text{if $s^{(j)}=1$,} \\
[f(\la^{(j)},\mu^{(j)}),g(\la^{(j)},\mu^{(j)})]^*,&\text{if
$s^{(j)}=-1$.}
\end{cases}
\end{equation*}
Figure 1 illustrates this procedure for a pair $(\lambda,\mu)$ of
tableaux of size 8.

\begin{center}
\epsfig{figure=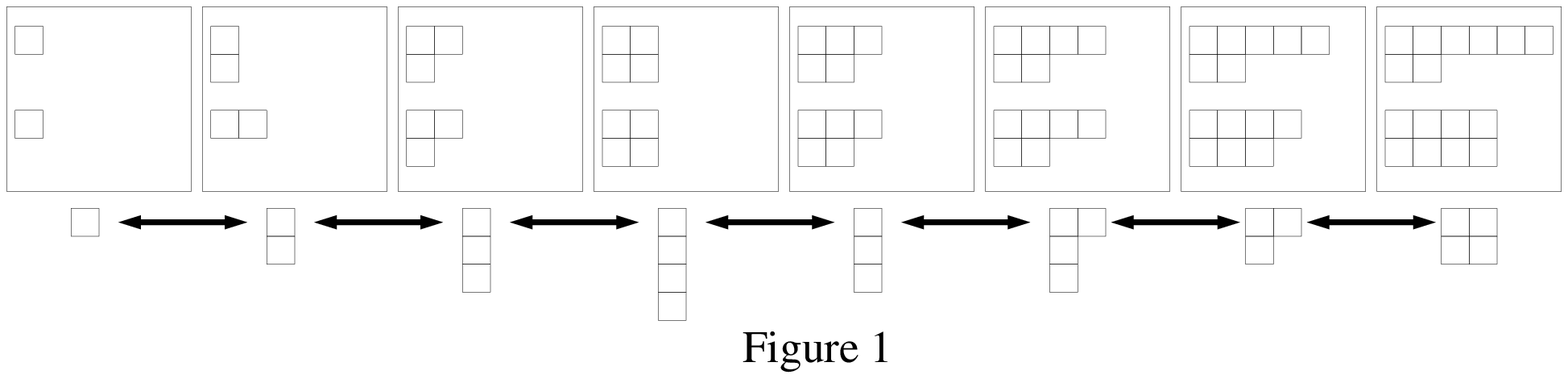,width=5in}
\end{center}

We need to check that
$$[0]=\nu^{(0)}\lra\nu^{(1)}\lra\cdots\lra\nu^{(m)}=\nu$$ with each $\nu^{(j)}\in\Ga(\ell)$.
Since $\Ga(\ell)$ is closed under $\ast$, to ensure that
$\nu^{(j)}\in\Ga(\ell)$ it is enough to show that
$[f(\la^{(j)},\mu^{(j)}),g(\la^{(j)},\mu^{(j)})]\in\Ga(\ell)$. This
holds because $\tilde \sigma_1+\tilde\sigma_2\le 4$ for every
diagram $\sigma$ with $\le 2$ rows, and
\begin{align*}
f(\la^{(j)},\mu^{(j)})+g(\la^{(j)},\mu^{(j)}) &= 2\max(\lambda_1^{(j)},\mu_1^{(j)}) - j \\
&= \max(\lambda_1^{(j)}-\lambda_2^{(j)},\mu_2^{(j)}-\mu_1^{(j)}) \le
\ell-2
\end{align*}
as $\la^{(j)},\mu^{(j)}\in\La(j,\ell)$. To show that
$\nu^{(j)}\lra\nu^{(j+1)}$ one checks the (four) cases corresponding
to the relationships $\la^{(j)}\rightarrow\la^{(j+1)}$ and
$\mu^{(j)}\rightarrow\mu^{(j+1)}$.  This is straightforward,
although tedious, with some care needed to see that the relationship
$\nu^{(j)}\lra\nu^{(j+1)}$ holds if $\nu^{(j)}$ has two rows and
$\nu^{(j+1)}$ has three.

Given an oscillating tableau $\nu^{(j)}$ in $\Gamma(\ell)$, we write
each $\nu^{(j)}$ as $[\nu_1^{(j)},\nu_2^{(j)}]$ or
$[\nu_1^{(j)},\nu_2^{(j)}]^*$. Let $s^{(j)}$ be $-1$ if and only if
there is a $\ast$ and $1$ if and only if there is not. This is
uniquely defined except when $\nu^{(j)}$ has exactly two rows. In
this case, we apply the following rule: if $k$ is the smallest
integer greater than or equal to $j$ such that $\nu^{(k)}$ has fewer
than $2$ or more than $2$ rows, then $s^{(j)} = s^{(k)}$; if there
is no such $k$, then $s^{(j)} = 1$. Figure 2 illustrates this
procedure:
\begin{center}
\epsfig{figure=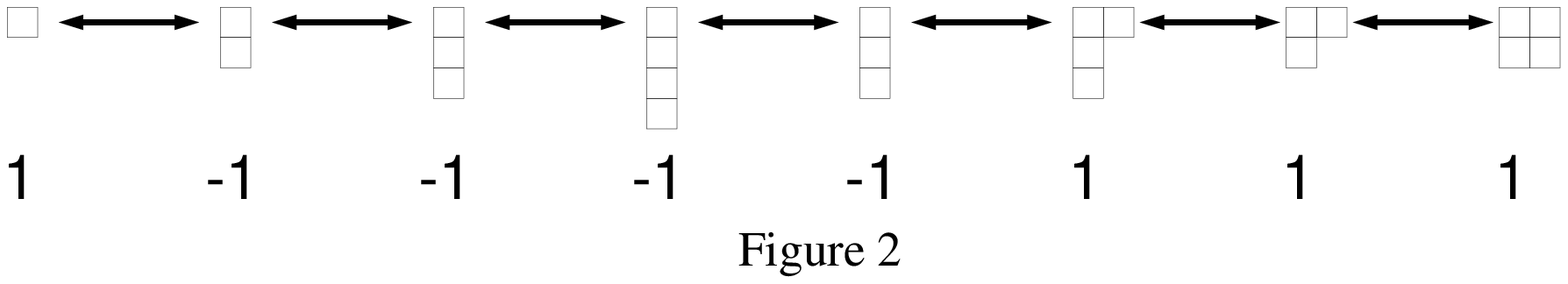,width=5in}
\end{center}

We next define for each $j$ between $1$ and $m$
\begin{align*}
\lambda_1^{(j)} &= \frac {j+\nu_1^{(j)}+s^{(j)}\nu_2^{(j)}}{2}, \\
\lambda_2^{(j)} &= j - \lambda_1^{(j)} \\
\mu_1^{(j)} &= \frac {j+\nu_1^{(j)}-s^{(j)}\nu_2^{(j)}}{2}, \\
\mu_2^{(j)} &= j - \mu_1^{(j)}.
\end{align*}
These are integers because
$$\nu_1^{(j)}+s^{(j)}\nu_2^{(j)}\equiv \nu_1^{(j)}+\nu_2^{(j)}\equiv |\nu^{(j)}|\equiv j\pmod2.$$
They are obviously all non-negative, and setting
$$\lambda^{(j)} = [\lambda_1^{(j)},\lambda_2^{(j)}],\ \mu^{(j)} = [\mu_1^{(j)},\mu_2^{(j)}],$$
we obtain diagrams in $\Lambda(j,\ell)$.  If $\nu^{(m)}$ has exactly
two rows, our sign convention guarantees that $s^{(m)}=1$ and
therefore that $\lambda_1 = \lambda_1^{(m)}> \mu_1^{(m)} = \mu_1$;
otherwise $\lambda_1=\mu_1$.

It is not difficult to see that these two constructions are mutually
inverse; the most delicate point is that signs $s^{(j)}$ are
respected.

The resulting bijection implies equation (\ref{tensordim})
immediately. Equations (\ref{symdim}) and (\ref{altdim}) follow by
setting $\lambda=\mu$ and defining an order on the tableaux
$t_\lambda$ of shape $\lambda$ according to which $t_\lambda\ge
t_\mu$ if and only if
$$(\lambda_1^{(m)},\lambda_1^{(m-1)},\ldots,\lambda_1^{(1)}) \ge
(\mu_1^{(m)},\mu_1^{(m-1)},\ldots,\mu_1^{(1)})$$
in lexicographic order.  Thus $t_\lambda\ge t_\mu$ if and only if
$\nu^{(n)}$ has $\le 1$ row and $t_\lambda < t_\mu$ if and only if
$\nu^{(n)}$ has $\ge 3$ rows.

\end{proof}

We remark that the above theorem allows us to deduce for each fixed
$\ell$ closed form expressions for $|\O_\ell(m,\la)|$ using the
corresponding expressions for $|\T_\ell(\la)|$. Such formulae can be
obtained via the interpretation of $|\T_\ell([m-p,p])|$ as the
dimension of a certain representations $V_{m,p}$ and
$\overline{V}_{m,p}$ (defined below) of the braid group $\B_m$.  For
example, we have

\begin{lemma}[Jones]
\begin{enumerate}

\item[(a)]For $\ell=\infty$,
$$|\T([m-p,p])|=\dim V_{m,p}=\binom{m}{p}-\binom{m}{p-1}$$
where $\binom{m}{-1}=0$ by convention.

\item[(b)]For $\ell=6$,
$$|\T_\ell([m-p,p])|=\dim(\overline{V}_{m,p})=\begin{cases} (3^{\lfloor\frac{m-1}{2}\rfloor}+1)/2 &
m-2p=0,1\\
(3^{\lfloor\frac{m-1}{2}\rfloor}-1)/2 & m-2p=3,4\\
3^{\lfloor\frac{m-1}{2}\rfloor} & m-2p=2\end{cases}$$
\end{enumerate}
\end{lemma}
\begin{proof}
The case with $\ell=\infty$ appears explicitly in \cite{J1}. The
$\ell=6$ formulae can be easily proved by induction using the
structure of $\Lambda(6)$---see the last example in
\cite[\S4.2]{J3}.
\end{proof}

The following technical lemma shows that pairs $(m,\la)$ with $\la\in\Ga(\ell)$ are
distinguished by pairs $(m-1,\nu)$ with $\nu\lra\la$ as long as $m\geq 3$, and will
be used in the proof of Theorem \ref{corresp}.  Define sets
$$P(m,\la):=\{\nu\in\Ga(\ell): \nu\lra\la,\quad (m-1)-|\nu|\in 2\N\}$$
of level $m-1$ predecessors of $\la$.  Then we have:
\begin{lemma}\label{diaglemma}
Fix $m\geq 3$, and let $\la, \mu\in\Ga(\ell)$ with
$m-|\la|,m-|\mu|\in 2\N$. Then $\la\neq\mu$ implies:
\begin{equation}\label{neighbors}
 P(m,\la)\neq P(m,\mu).\end{equation}

\end{lemma}
\begin{proof}
First suppose $|\la|\geq 3$, $|\mu|\geq 3$, $\nu\rightarrow\la$
implies $\nu\in\Ga(\ell)$ and $\nu\rightarrow\mu$ implies
$\nu\in\Ga(\ell)$. Then a direct application of \cite[Lemma
2.11(b)]{wenzl88} shows that (\ref{neighbors}) holds for all $\la,\mu$
with these restrictions.  The only diagrams that do not satisfy
these extra hypotheses are $[0],[1],[2],[1,1]$ and $[\ell-2,1,1]$
(which fails because $[\ell-2,1]\not\in\Ga(\ell)$). The diagram
$[\ell-3,1,1]$ is the unique diagram in $\Ga(\ell)$ such that
$[\ell-2,1,1]\lra[\ell-3,1,1]$, so clearly (\ref{neighbors}) holds for
$\la=[\ell-3,1,1]$ and $\mu$ arbitrary.  The remaining diagrams can
be handled similarly, noting that since $m\geq 3$,
$[1,1]\lra[1,1,1]$, $[2]\lra[3]$, $[0]\lra [1]$ and $[1]$ is the
unique diagram with the latter property.
\end{proof}
\section{Temperley-Lieb Algebras}

Temperley-Lieb algebras are natural representation spaces for braid
groups. They admit a natural trace.  The Jones polynomial of a link
$L$ is defined as the trace of any braid $\beta$ for which the
corresponding closed braid $\hat\beta = L$.

Fix a complex variable $q$.
\begin{defn}
The Temperley-Lieb algebra $T_m(q)$
 is the $\C(q)$-algebra generated by $e_1,e_2,\cdots,e_{m-1}$ satisfying:
 \begin{enumerate}
 \item[(T1)] $e_ie_{i\pm 1}e_i=\frac{q^{-2}}{(1+q^{-2})^2}e_i=\frac{1}{(q+q^{-1})^2}e_i$
 \item[(T2)] $e_ie_j=e_je_i$ for $|i-j|\geq 2$
 \item[(H)] $e_i^2=e_i$
 \end{enumerate}\end{defn}
 By convention we put $T_0=T_1=\C(q)$.
When there is no danger of confusion we will denote $T_m(q)$ simply
by $T_m$.
 \begin{remark}
 The reader is warned that our definition of the Temperley-Lieb algebra differs slightly from
 the standard one (see \cite{GWTL} for example) in
 which the Temperley-Lieb algebras are defined with parameter
 $t$ which corresponds to $q^{-2}$ in our definition.
 \end{remark}

 \begin{lemma}\label{tlrels} Define $g_i:=(1+q^{-2})e_i-1$.
 \begin{enumerate}
 \item[(a)] The following relations hold in $T_m$:
 \begin{enumerate}
 \item[(B1)] $g_ig_{i+1}g_i=g_{i+1}g_ig_{i+1}$ for $1\leq i\leq m-2$
\item[(B2)] $g_ig_j=g_jg_i$ if $|i-j|\geq 2$
 \item[(T3)] $g_i^{-1}=(q^{2}+1)e_i-1$
\item[(T4)] $g_ie_i=q^{-2}e_i$
\item[(T5)] $e_ig_{i+1}e_i=\frac{-1}{q^{-2}+1}e_i$
\item[(T6)] $g_ig_{i\pm 1}g_i+g_ig_{i\pm 1}+g_{i\pm 1}g_i+g_i+g_{i\pm 1}+1=0$
\item[(T7)] $(g_i+1)(g_i-q^{-2})=0$.
 \end{enumerate}
 \item[(b)]
 The inductive limit of the algebras $T_m$ admits a $\C(q)$-valued trace $\tr$ uniquely determined by:
 \begin{enumerate}
 \item[(M1)] $\tr(1)=1$
 \item[(M2)] $\tr(ab)=\tr(ba)$
 \item[(M3)] $\tr(ae_{m-1})=\frac{1}{(q+q^{-1})^2}\tr(a)$ for $a\in T_{m-1}$.
 \end{enumerate}
 \end{enumerate}
 \end{lemma}
The relations (B1) and (B2) imply that $T_m$ is a quotient of
$\C(q)\B_m$.  Moreover, one deduces from these relations that $T_m$
is finite-dimensional over $\C(q)$.

Specializations of Temperley-Lieb algebras remain well-defined for
$q\not\in\{0,\pm i\}$ from which we obtain $\C$-algebras and
$\C$-representations of
 $\B_m$ factoring over $T_m$.  The analysis of these specializations
 breaks naturally into two cases: 1) the generic case--those $q$ for
 which $q^{2k}-1\in\C^*$ for all integers $k\geq 1$ and 2) the
 proper root of unity case--those $q$ for which $q^2$ is a primitive
 $\ell$th root of unity with $\ell\geq 3$.  When we wish to consider
 both cases simultaneously we say that $q^2$ is a primitive $\ell$th
 root of unity with $3\leq\ell\leq\infty$ where the case
 $\ell=\infty$ covers that former case.  By an abuse of notation we
 will continue to denote these specializations by $T_m$ since $\ell$
 and $q$ will always be clear from the context.

In the generic case the the trace $\tr$ is faithful, \emph{i.e.} the
annihilator ideal $J_m:=\{a\in T_m: \text{$\tr(ab)=0$ for all $b\in
T_m$}\}=\{0\}$.  Moreover, in these cases the algebras $T_m$ are
semisimple.

When $q^2$ is a primitive $\ell$th root of unity for
$3\leq\ell<\infty$ the specializations are not semisimple, and $\tr$
is not faithful.  However, the annihilator of the trace
$J_{m}(q,\ell)$ contains the Jacobson radical and the (semisimple)
quotient algebra $T_m/J_{m}(q,\ell)$
 will be denoted by $\overline{T}_m$.

   As semisimple finite
dimensional algebras, $T_m$ and $\TL$ are direct sums of full matrix
algebras.  The simple subalgebras of $T_m$ and $\TL$ are in
one-to-one correspondence with the subsets $\La(m,\ell)\subset YD$,
where $\ell=\infty$ covers the generic case. The decompositions of
$T_m$ and $\TL$ into full matrix algebras and the restriction rules
are described in the following:
\begin{prop}
\label{tldecomp} Define $T_{2,0}$ and $T_{2,1}$ to be the
eigenspaces of $g_1\in T_2$ corresponding to eigenvalues $-1$ and
$q^{-2}$ respectively.
\begin{enumerate}
\item[(a1)] For the generic cases $\ell=\infty$, we have: $$T_m=\bigoplus_p T_{m,p}$$ where
$0\leq p\leq\lfloor\frac{m}{2}\rfloor$, and $T_{m,p}$ is a full
matrix algebra, corresponding to Young diagram $[m-p,p]$.
\item[(a2)] For any Young diagram $[m-p,p]$, denote by $V_{m,p}$ the
$T_m$-representation such that $T_{m,p}\cong\End(V_{m,p})$.  Then
the restriction of $V_{m,p}$ to $T_{m-1}$ decomposes irreducibly as
$$V_{m-1,p}\oplus V_{m-1,p-1}$$ where we
set $V_{n,t}=\{0\}$ if $[n-t,t]\not\in\La(n,\infty)$.
\item[(b1)] Suppose $q^2$ is an $\ell$ root of unity with $3\leq\ell<\infty$.
Then $$\TL=\bigoplus_{p}\overline{T}_{m,p}$$ where the sum is over
all $p$ such that $[m-p,p]\in\Lambda(m,\ell)$ and
$\overline{T}_{m,p}(\ell)$ is a full matrix algebra.
\item[(b2)] Denote by $\overline{V}_{m,p}$ the
$\TL$-representation such that
$\overline{T}_{m,p}\cong\End(\overline{V}_{m,p})$ and set
$\overline{V}_{n,t}=\{0\}$ if $[n-t,t]\not\in\Lambda(n,\ell)$. Then
the restriction of $\overline{V}_{m,p}$ to $\overline{T}_{m-1}$
decomposes irreducibly as:
$$\overline{V}_{m-1,p}\oplus \overline{V}_{m-1,p-1}$$ where we discard any
summand that is $\{0\}$.
\end{enumerate}
\end{prop}

The restriction rule given above can be more easily explained
combinatorially: the representation $V_{m-1,s}$ appears in the
restriction of $V_{m,p}$ to $T_{m-1}$ if and only if
$[m-1-s,s]\rightarrow[m-p,p]$ with an analogous statement for $\TL$.
 From this description we obtain the \emph{Bratteli diagrams} for
$T_m$ and $\TL$.  These are graphs with vertices labelled by
diagrams $[m-p,p]\in\La(\ell)$ (where $\ell=\infty$ covers the
generic case as usual) and with an edge between the vertex labelled
by $[m-p,p]$ and $[m-1-s,s]$ if and only if $V_{m-1,s}$ is a
$T_{m-1}$-subrepresentation of $V_{m,p}$ with the same statement for
$\TL$ and $\overline{V}_{m,p}$.  The ambiguity with the two
components of $T_2$ is removed by the definition of $T_{2,0}$ and
$T_{2,1}$ as above.  Observing that $\dim(T_0)=\dim(T_1)=1$, we can
inductively compute the dimensions of $V_{m,p}$ and
$\overline{V}_{m,p}$ by counting the increasing paths $t_{[m-p,p]}$
in the Bratteli diagram of $T_m$ or $\TL$ from $[0]$ to $[m-p,p]$.
Thus the representation spaces $V_{m,p}$ and $\overline{V}_{m,p}$
have bases labelled by such paths:
$$[0]\rightarrow[1]\rightarrow\cdots\rightarrow[m-p,p],$$ where each diagram must be in $\Lambda(\ell)$.
But increasing paths in the Bratteli diagrams of $T_m$ and $\TL$ are
just Young tableaux restricted to $\La(\ell)$ (for $\ell=\infty$ and
$6\leq\ell<\infty$ respectively) so from this we see that
$$\dim V_{m,p}=|\T_{\infty}([m-p,p])|; \dim \overline{V}_{m,p}=|\T_{\ell}([m-p,p])|.$$

\section{BMW Algebras}
While the Jones polynomial $V_L(t)$ was derived from the trace on
the Temperley-Lieb algebras, the two-variable Kauffmann polynomial
$F_L(a,z)$ \cite{kauff} was first defined in a purely combinatorial
way. However, not long after its definition, Birman-Wenzl and
Murakami (\cite{BW}, \cite{mur}) independently found the appropriate
traced quotients of the braid group algebras corresponding to
$F_L(a,z)$, and they are now known as BMW (or $q$-Brauer) algebras.
The reader is warned that the parameters $r$ and $q$ below
correspond to a different version $K_L(r,q)$ of the Kauffmann
polynomial related to $F_L(a,z)$ by a non-trivial change of
variables.
\subsection{Definitions and Algebraic Results}
\begin{defn}
The BMW algebra $\CC_m(r,q)$ is the $\C(r,q)$-algebra with
invertible generators $G_1, G_2,\ldots, G_{m-1}$ satisfying the
braid relations $(B1)$ and $(B2)$ above and:
\begin{enumerate}
\item[(R1)] $(G_i-r^{-1})(G_i-q)(G_i+q^{-1})=0$
\item[(R2)] $E_iG_{i-1}^{\pm 1}E_i=r^{\pm 1}E_i$, where
\item[(E)] $(q-q^{-1})(1-E_i)=G_i-G_i^{-1}$ defines $E_i$.
\end{enumerate}
\end{defn}
By convention $\CC_0(r,q)=\CC_1(r,q)=\C(r,q)$. These relations
imply:

\begin{prop}\cite[\S 3]{BW}\label{BMWspan}
\begin{enumerate}
\item[(a)]
The algebra $\CC_m(r,q)$ is linearly spanned by elements of the form
$a\chi b$ where $a,b\in\CC_{m-1}(r,q)$ are monomials and
$\chi\in\{1,G_{m-1},E_{m-1}\}$.
\item[(b)] The elements in $\CC_m(r,q)$ spanned by monomials of the form $aE_{m-1}b$ with
$a,b\in\CC_{m-1}(r,q)$ form an ideal $\I_m$, and $\CC_m(r,q)/\I_m$
is isomorphic to the Hecke algebra $\Hk_m(q^2)$.
\end{enumerate}
\end{prop}

The inductive limit of the algebras $\CC_m(r,q)$ is equipped with a
 trace:

\begin{prop}\cite[Lemma 3.4]{We}\label{BMWtr}
Set $x=\frac{r-r^{-1}}{q-q^{-1}}+1$. There exists a functional,
$\Tr$, on $\CC_\infty(r,q)$ uniquely defined inductively by:
\begin{enumerate}
\item $\Tr(1)=1$
\item $\Tr(ab)=\Tr(ba)$
\item $\Tr(E_i)=1/x$
\item $\Tr(G_i^{\pm 1})=r^{\pm 1}/x$
\item $\Tr(a\chi b)=\Tr(\chi)\Tr(ab)$ for $a,b\in\CC_{m-1}(r,q)$, $\chi\in\{G_{m-1},E_{m-1}\}$.
\end{enumerate}
\end{prop}
As in the case of Temperley-Lieb algebras, one may specialize $r$
and $q$ to be complex numbers and for any specialization for which
$\CC_m(r,q)$ and $\Tr$ are well-defined both Propositions
\ref{BMWspan} and \ref{BMWtr} still hold. For such $r$ and $q$
denote the annihilator ideal of $\Tr$ on $\CC_m(r,q)$ by:
$A_m(r,q)$. As long as $r\not=\pm q^n$ for any integer $n$ and
$q^{2k}-1\in\C^*$ for all $k\geq 1$, the trace $\Tr$ is faithful on
$\CC_m(r,q)$; moreover, $\CC_m(r,q)$ is semisimple.  We will shorten
$\CC_m(r,q)$ to $\CC_m$ for these generic cases.

The specializations of BMW algebras with $r=\pm q^n$ are related
(via quantum Schur-Weyl-Brauer duality) to quantum groups of Lie
types $B, C$ and $D$, while if we further specialize $q^2$ to be a
primitive $\ell$th root of unity we obtain interesting
$\C$-representations of $\B_m$ in analogy with the Temperley-Lieb
situation. When $r=q^n$ and/or $q^2$ is a root of unity, the BMW
algebras fail to be semisimple.  By taking the quotient by the ideal
$A_m(r,q)$ semisimplicity can often be recovered.  For example,
\begin{prop}[\cite{We}]
Fix $r$ and $q$ with $r=q^n$ where $3\leq n\leq\ell-3$ and $q^2$ is
a primitive $\ell$th root of unity with $\ell\leq\infty$.  Then
$\overline{\CC}_m(r,q):=\CC_m(r,q)/A_m(r,q)$ is semisimple.
\end{prop}
As usual, we designate the case where $q^{2k}-1\in\C^*$ for all
$k\geq 1$ by $\ell=\infty$.

\subsection{Representation Theory}
In this paper we are interested in the cases where $r=q^3$ and $q^2$
is a primitive $\ell$th root of unity with $\ell\leq\infty$. As
described in \cite{LRW} Prop. 6.2~(1)(c), the
 $\B_m$-representations factoring over $\CC_m(q^3,q)$ are
non-degenerate provided $q^2$ is an $\ell$th root of unity with
$6\leq\ell\leq\infty$. The simple subalgebras of $\CC_m$ are in
one-to-one correspondence with Young diagrams $\la$ with $m-|\la|\in
2\N$, while for the semisimple quotients of the specializations of
$\overline{\CC}_m(r,q)$ with $r=q^3$ and $q^2$ an $\ell$th root of
unity with $6\leq\ell\leq\infty$ one must restrict to diagrams in
the set $\Gamma(\ell)$
 defined above.

We have the following description of the simple decompositions and
restriction rules for BMW algebras in both the generic case and
the specializations we study:
\begin{prop}\label{bmwrep}Define $\CC_{2,[0]}$, $\CC_{2,[1,1]}$  and $\CC_{2,[2]}$ to be the
eigenspaces of $G_1\in \CC_2$ corresponding to eigenvalues $r^{-1}$,
$-q^{-1}$ and $q$ respectively.
\begin{enumerate}
\item[(a1)] Suppose $r\not=\pm q^n$ and $q^2$ is not a root of unity.
Then $\CC_m=\bigoplus_{\la} \CC_{m,\la}$ where $m-|\la|\in 2\N$, and
$\CC_{m,\la}$ is a full matrix algebra.
\item[(a2)] If $W_{m,\la}$ is a simple $\CC_{m,\la}$-module then the
restriction of $W_{m,\la}$ to $\CC_{m-1}(r,q)$ decomposes
irreducibly as:
$$\bigoplus_{\mu\lra\la}W_{m-1,\mu}.$$
\item[(b1)] Suppose $r=q^3$ and $q^2$ is an $\ell$th root of unity with $6\leq\ell\leq\infty$.
Then
$$\overline{\CC}_m(r,q)=\bigoplus_{\la}\overline{\CC}_{m,\la}(r,q)$$
where $m-|\la|\in 2\N$ and $\la\in\Gamma(\ell)$ and
$\overline{\CC}_{m,\la}(r,q)$ is a full matrix algebra.
\item[(b2)] Let $\overline{W}_{m,\mu}$ be a simple $\overline{\CC}_m(r,q)$ module
with $r$ and $q$ as in \textrm{(b1)}. Then the restriction of
$\overline{W}_{m,\mu}$ to $\overline{\CC}_{m-1}(r,q)$ decomposes
irreducibly as:
$$\bigoplus_{\stackrel{\mu\lra\la}{\mu\in\Gamma(\ell)}}\overline{W}_{m-1,\mu}.$$
\end{enumerate}
\end{prop}

This description gives us a convenient way of encoding the
inclusions of BMW algebras via their Bratteli diagrams.  The
ambiguity between the three simple components for $m=2$ is removed
by assigning the labels to eigenspaces as in the proposition above.
Define a graph whose vertices are labelled by $(m,\la)$ where
$m-|\la|\in 2\N$ and the labels $(m,\la)$ and $(m-1,\mu)$ are
connected by an edge if and only if $\la\lra\mu$.  For
specializations of $\overline{\CC}_m(r,q)$ with $r=q^3$ and $q^2$
and $\ell$th root of unity with $6\leq\ell\leq\infty$, the Bratteli
diagram is defined in the same way except that the Young diagrams
are restricted to be in the set $\Ga(\ell)$ defined in Section
\ref{comb}.  From this we see that there are bases for $W_{m,\la}$
and $\overline{W}_{m,\la}$ indexed by the set of paths of length $m$
in the Bratteli diagram beginning at $[0]$ and ending at $\la$
(where all diagrams must be in $\Ga(\ell)$ in the latter case). From
the structure of the Bratteli diagrams we see that these paths are
in one-to-one correspondence with oscillating tableaux. Thus, the
dimension of $W_{m,\la}$ (respectively, $\overline{W}_{m,\la}$) is
the number $|\O(m,\la)|$ of oscillating tableaux (resp.
$|\O_{\ell}(m,\la)|$) of shape $\la$ and length $m$. Note that when
$r=q^3$ and $q^2$ is a primitive $\ell$th root of unity the Bratteli
diagram for $\overline{\CC}_m(r,q)\subset\overline{\CC}_{m+1}(r,q)$
depends only on $\ell$, not on the specific choice of $q$.

\section{Symmetric Squares of Algebras}

Let $A$ be an associative $\C$-algebra.   We define $\s^2 A$ to be
the subalgebra of $A\otimes_{\C} A$ generated by $\{a\otimes a\mid
a\in A\}$.

\begin{lemma}
Let $A = A_1\oplus\cdots \oplus A_n$.  Then
$$\s^2 A = \bigoplus_{i=1}^n \s^2 A_i \oplus \bigoplus_{j=1}^{n-1}\bigoplus_{k=j+1}^n
A_j\otimes A_k.$$
\end{lemma}

\begin{proof}
Let $\sigma_{jk}\colon A_j\otimes A_k\to A_k\otimes A_j$ exchange
factors. The natural inclusion
$$\iota\colon (x_i,y_{jk})\mapsto \sum_{i=1}^n x_i + \sum_{j=1}^{n-1}\sum_{k=j+1}^n (y_{jk}+\sigma_{jk}(y_{jk}))$$
is obviously injective.  It is surjective because for all
$a=a_1+\cdots+a_n$, with $a_i\in A_i$, we have
$$
a\otimes a = \iota(a_1\otimes a_1,\ldots,a_n\otimes a_n,a_1\otimes
a_2,\ldots,a_{n-1}\otimes a_n).
$$
\end{proof}

\begin{prop}
\label{blocks} If $A=M_n(\C)$, then
$$\s^2 A = M_{\binom {n+1}2}(\C)\oplus M_{\binom{n}{2}}(\C).$$
If $B=M_m(\C)$, then $A\otimes B = M_{mn(\C)}$.
\end{prop}

\begin{proof}
The proposition is trivial when $n=1$ (where $M_0(\C)$ is understood
to mean the zero-ring). We therefore assume $n\ge 2$.

If $V = \C^n$, then $\GL_n(\C)$ acts on $V$, and $V\otimes V$
decomposes as a direct sum of two irreducible
$\GL_n(\C)$-representations: $\s^2 V$ and $\wedge^2 V$.   Thus, the
diagonal image of $\GL_n(\C)$ in $\End(V\otimes V)$ lies in
$$\End(\s^2 V)\oplus \End(\wedge^2 V)\cong M_{\binom {n+1}2}(\C)\oplus M_{\binom{n}{2}}(\C).$$
As $\GL_n(\C)$ is dense in $M_n(\C)$, the same is true of the
diagonal image of $M_n(\C)$, and it follows that the subalgebra of
$\End(V\otimes V)$ generated by $a\otimes a$ for $a\in M_n(\C)$ is
contained in $\End(\s^2 V)\oplus \End(\wedge^2 V)$.  Conversely, if
$\s^2 M_n(\C)$ is a $\ast$-subalgebra of $M_{n^2}(\C)$ and therefore
a semisimple algebra.  If it is properly contained in $\End(\s^2
V)\oplus \End(\wedge^2 V)$, then it has a larger centralizer in
$\End(V\otimes V)$, so the centralizer of the diagonal image of
$\GL_n(\C)$ in $\End(V\otimes V)$ has dimension $>2$.  This is
impossible by Schur's lemma; we have already observed that $V\otimes
V$ decomposes as the sum of two inequivalent irreducible
representations of $\GL_n(\C)$.

For the second claim, let $W = \C^m$.   There is a natural map
$\End(V)\otimes \End(W)\to \End(V\otimes W)$ which is an isomorphism
since
\begin{multline*}
\End(V)\otimes \End(W) =
(V\otimes V^*) \otimes (W\otimes W^*) \\
= V\otimes W\otimes V^*\otimes W^* = \End(V\otimes W).
\end{multline*}
\end{proof}

\begin{remark}
A natural setting in which to consider symmetric squares of algebras
is that of von Neumann algebras.  The second part of
Proposition~\ref{blocks} is well known to hold for factors.  We do
not know whether the first part holds as well, i.e., whether the
symmetric square of a non-trivial factor is always the direct sum of
two factors.
\end{remark}

If $A$ is endowed with a linear functional $\tr\colon A\to \C$
satisfying the trace identity $\tr(ab) = \tr(ba)$, then
$\tr\otimes\tr\colon A\otimes A\to \C$ also satisfies the trace
identity, so the same is true of its restriction (denoted $\tr^2$)
to $\s^2 A$.

We apply the symmetric square construction to Temperley-Lieb
algebras. As $T_m=\bigoplus_i T_{m,i}$ is a direct sum of full
matrix algebras, the symmetric square $\s^2T_m$ can be decomposed
as:
$$\bigoplus_i \s^2T_{m,i}\oplus \bigoplus_{j<k} T_{m,j}\otimes T_{m,k}.$$
Thus $V_{m,j}\otimes V_{m,k}$ ($j<k$), $\s^2V_{m,i}$ and
$\bigwedge^2 V_{m,j}$ are irreducible representations of $\ST$.  The
trace $\tr$ on $T_m$ determines the trace $\tr^2$ on $\ST$.

Observe that the same analysis applies to the symmetric square
$\s^2\TL$ for $\ell<\infty$ with analogous conclusions replacing
$V_{m,i}$ by $\overline{V}_{m,i}$ and restricting to
$\Lambda(\ell)$-diagrams in all formulae.

Fix $q$ such that $q^2$ is an $\ell$th root of unity with
$6\leq\ell\leq\infty$ and set $r=q^3$ and
$x=\frac{r-r^{-1}}{q-q^{-1}}+1=(q+q^{-1})^2$.  By an abuse of
notation we will continue to denote the images of the generators of
$T_m$ in $\TL$ by $g_i$ and $e_i$ as this should cause no confusion.
We define elements $\tilde{G}_i=q(g_i\otimes g_i)$ and
$\tilde{E}_i:=x(e_i\otimes e_i)$
 of $\ST$ or $\s^2\TL$ and derive
some relations from those of $T_m$:
\begin{lemma}\label{relslemma}We have the following identities:
\begin{enumerate}
\item $(\tG-r^{-1})(\tG-q)(\tG+q^{-1})=0$
\item $(q-q^{-1})(1-\tE)=(\tG-\tG^{-1})$
\item $\tE_i\tG^{\pm 1}_{i-1}\tE_i=r^{\pm 1}\tE_i$
\item $\tr^2(1)=1$
\item $\tr^2(ab)=\tr^2(ba)$
\item $\tr^2(\tE_i)=1/x$
\item $\tr^2(\tG_i)=r^{\pm 1}/x$
\item $\tr^2(a\chi b)=\tr^2(\chi)tr^2(ab)$ for $a,b\in\ST$,
$\chi\in\{\tG_m,\tE_m\}$
\end{enumerate}
\end{lemma}
\begin{proof}All of these relations follow directly from Lemma
\ref{tlrels}.  For example, let us prove (2).
\begin{align*}\tG_i-\tG_i^{-1}&=q(g_i\otimes g_i)-(g_i^{-1}\otimes
g_i^{-1})/q\\
&=q((q^{-2}+1)^2 e_i\otimes e_i-(q^{-2}+1)(1\otimes e_i+e_i\otimes 1)+1\otimes 1)\\
&\quad-q^{-1}((q^{2}+1)^2 e_i\otimes e_i-(q^{2}+1)(1\otimes
e_i+e_i\otimes 1)+1\otimes 1)\\
&=(q(q^{-2}+1)^2-(q^2+1)^2/q)e_i\otimes e_i+(q-q^{-1})1\otimes 1\\
&=(q-q^{-1})(1\otimes 1-(q+q^{-1})^2 e_i\otimes e_i)\\
&=(q-q^{-1})(1-\tE_i)
\end{align*}
\end{proof}

\section{The Isomorphism}
Throughout this section set $r=q^3$ and let $q^2$ be an $\ell$th
root of unity with $6\leq\ell\leq\infty$. Consider the mapping:
$$\Phi(G_i)=\tG_i,\quad 1\leq i\leq m-1$$  It is immediate from Lemma \ref{relslemma}
and the defining relations of $\CC_m(r,q)$ that $\Phi$ extends to an
algebra homomorphism $\CC_m(r,q)\rightarrow\ST$.
 Another consequence of Lemma \ref{relslemma}
is that $\Phi(E_i)=\tilde{E}_i$.   We can now prove:
\begin{lemma}\label{injlemma}   The induced map
$$\overline{\Phi}:\overline{\CC}_m(r,q)=\CC_m(r,q)/A_m(r,q)\rightarrow
\Phi(\CC_m(r,q))/\Phi(A_m(r,q))$$ is injective.
\end{lemma}
\begin{proof}
It is enough to show that $\ker \Phi\subset A_m(r,q)$. First note
that $\tr^2$ induces a trace form $\Phi^{-1}(\tr^2)$ on $\CC_m(r,q)$
that has the Markov property and the values of $\Phi^{-1}(\tr^2)$
and $\Tr$ coincide on $\{1,E_i,G_i\}$ for all $i$ so that the
uniqueness of $\Tr$ implies that $\Phi^{-1}(\tr^2)=\Tr$.  Suppose
$a\in \ker\Phi$, and $b\in\CC_m(r,q)$.  Then
$\Tr(ab)=\tr^2(\Phi(ab))=\tr^2(0)=0$ so that $a\in A_m(r,q)$.
\end{proof}

An immediate corollary of this lemma is that $\overline{\CC}_m(r,q)$
is isomorphic to a semisimple quotient of the subalgebra of $\ST$ or
$\s^2\TL$ generated by $\{\tG_i\}$ so that
$\dim(\overline{\CC}_m(r,q))\leq \dim(\ST)$ for $\ell=\infty$ and
$\dim(\overline{\CC}_m(r,q))\leq\dim(\s^2\TL)$ for
$6\leq\ell<\infty$.  But by Theorem \ref{maincomb},
\begin{eqnarray*}
\dim(\overline{\CC}_m(r,q))&=&\sum_{\stackrel{\nu\in\Ga(\ell)}{m-|\nu|\in 2\N}}|\O_\ell(m,\nu)|^2\\
&\geq&\sum_{p<q\le \frac m2}|\T_\ell([m-p,p])|^2\cdot|\T_\ell([m-r,r])|^2\\
&+&\sum_{p\le \frac m2}
\binom{1+|\T_\ell([m-p,p])|}{2}^2\\
&+&\sum_{p\le \frac m2}
\binom{|\T_\ell([m-p,p])|}{2}^2\\
&=&\begin{cases}\dim(\ST) & \ell=\infty\\
\dim(\s^2\TL) & 6\leq\ell<\infty, \end{cases}\end{eqnarray*} so by
dimension we have our main result:
 \begin{theorem} Let $r=q^3$.  Then
 $\overline{\CC}_m(r,q)\cong\ST$ if $q$ is not a root
 of unity and
$\overline{\CC}_m(r,q)\cong\s^2\TL$ for $q^2$ an $\ell$th root of
unity with $6\leq\ell<\infty$.
\end{theorem}

Although we have established isomorphisms between these semisimple
algebras as promised, we have not identified the images of the
simple components of $\overline{\CC}_m(r,q)$ under
$\overline{\Phi}$. Not surprisingly, the combinatorial
correspondence in Theorem \ref{maincomb} is compatible with
$\overline{\Phi}$:
\begin{theorem}\label{corresp}
Fix $[m-s,s],[m-t,t]\in\La(m,\ell)$ with and define $\nu_1$ and
$\nu_2$ as in Theorem \ref{maincomb}.  Then the map
$\overline{\Phi}$ induces isomorphisms of simple algebras as follows
for $6\leq\ell<\infty$:
\begin{enumerate}
\item if $\la\neq\mu$, $\overline{\CC}_{m,[\nu_1,\nu_2]}(r,q)\cong\End(\oV_{m,s})\otimes\End(\oV_{m,t})$
\item if $\la=\mu$, $\overline{\CC}_{m,[\nu_1]}(r,q)\cong \End(\s^2\oV_{m,s})$,
\item if $\la=\mu$, $\overline{\CC}_{m,[\nu_1]^*}(r,q)\cong \End(\wedge^2\oV_{m,s})$
\end{enumerate}
The same statement holds for $\ell=\infty$ replacing $\overline{\CC}_{m,\lambda}$
and $\overline{V}$ by $\CC_{m,\lambda}$ and $V$ respectively.
\end{theorem}
\begin{proof}The cases $m=0,1$ are clear since all algebras in
question are isomorphic to $\C$.  For $m\geq 2$ we proceed by
induction on $m$.  The (base) case $m=2$ follows by checking that
the labelling conventions for the eigenspaces of $\tG_1\in\s^2T_2$ (induced from
those of $g_1\in T_2$) and $G_1\in\CC_2$ are compatible
with the correspondence of Theorem \ref{maincomb}.  
Now suppose that
the statement holds for some $m-1\geq 2$. By Theorem \ref{bmwrep}
and Lemma \ref{diaglemma} any simple component of
$\overline{\CC}_m(r,q)$ is determined by the set of labels 
of the simple $\overline{\CC}_{m-1}(r,q)$-subalgebras
contained in $\overline{\CC}_m(r,q)$.  Applying the induction
hypothesis to these simple $\overline{\CC}_{m-1}(r,q)$-subalgebras
we obtain isomorphisms between the simple components of
$\s^2\overline{T}_{m-1}$ and those of $\overline{\CC}_{m-1}(r,q)$ as
in the statement of the theorem.  Tracing through the corresponding
labels of the simple components we see that this implies the result
for $m$.
\end{proof}

\section{Braid Group Images}

The irreducible representations of $\B_m$ factoring over $\TL$ are
unitary if $q=e^{\pm\pi i/\ell}$. The closed images of these unitary
$\B_m$-representations have been classified in \cite{J1}, \cite{BWj}
and \cite{FLW}.  Our goal in this section is to solve analogous
problem for BMW algebras when $r=q^3$.  This was the original
motivation of this paper. The question of unitarity for
representations of $\B_m$ factoring over $\CC_m(r,q)$ is not so
simple in general.  It was shown in \cite{R} that the cases $r=q^n$
with $n<0$ even and $q$ \emph{any} primitive $\ell$th root of unity
with $\ell$ odd can fail to yield unitary representations of $\B_m$.
However, Wenzl \cite{We} showed that for essentially all other
$r=q^n$ with $q=e^{\pm\pi i/\ell}$ one obtains unitary $\B_m$
representations; in particular, this is so for $r=q^3$ and $q=e^{\pm\pi
i/\ell}$ with $6\leq\ell$.

Throughout this section, we will fix an integer $\ell\ge 6$ and
assume $q=e^{\pm \pi i/\ell}$ and $r=q^3$.  By Proposition~\ref{bmwrep},
we have a decomposition
$$\CC_m(r,q) = \bigoplus_{\lambda\in\Gamma(\ell)}\End(\overline{W}_{m,\lambda}).$$
Let $\sigma_{m,\lambda}\colon \B_m\to \bar \GL(\overline{W}_{m,\lambda})$ denote the corresponding
representation and $\bar\sigma_{m,\lambda} \colon \B_m\to \bar \PGL(\overline{W}_{m,\lambda})$
its projectivization.
Because of our choice of $(q,r)$, $\sigma_{m,\lambda}$ is always unitary.
The topological closure $\bar H_{m,\lambda}$ of
$\bar\sigma_{m,\lambda}(\B_m)$ is therefore a compact Lie group.

We recall that by Proposition~ \ref{tldecomp},
$$\overline{T}_m = \bigoplus_s \End(\overline{V}_{m,s}),$$
where $(\rho_{m,s},\overline{V}_{m,s})$ is the Jones representation corresponding to the Young diagram $[m-s,s]$.
Let $\bar G_{m,s}$ denote the closure of the image of the projectivized Jones representation
$\bar\rho_{m,s}(\B_m)$.  As $\rho_{m,s}$ is unitary, this is a compact Lie group.

We recall the precise result:

\begin{prop}\label{tlimages}
Let $m\ge 3$ and $\ell\ge 6$ be integers, and let $q=e^{\pm \pi i/\ell}$.  Let
$1\le s\le m/2$, and set
$d_{m,s}:=\dim \oV_{m,s}$.
Then
\begin{enumerate}
\item[(a1)]\cite{BWj,FLW}
If $\ell=6$ and $m$ is odd, then $d_{m,s}=\frac{3^{(m-1)/2}\pm1}{2}$ and
$\bar G_{m,s}\cong\PSp_{m-1}(3)$.
\item[(a2)]\cite{BWj} If $\ell=6$ and $m$ is even, then $d_{m,s}\in\{\frac{3^{(m-2)/2}\pm
1}{2},3^\frac{m-2}{2}\}$ and
$$\bar G_{m,s}\cong\begin{cases}\PSp_{m-2}(3)\ltimes (\Z_3)^{m-2}& s = m/2-1\\
\PSp_{m-2}(\Z_3)& s\in\{m/2-2,m/2\}.\end{cases}$$
\item[(b)]\cite{J1} If $\ell=10$ then $\bar G_{3,1}\cong \bar G_{4,2}\cong A_5$.
\item[(c)]\cite{FLW} Except in cases (a) and (b), $\bar G_{m,s}=\PSU(d_{m,s})$.
\item[(d)]\cite{FLW} If $s\neq t$, then the $\B_m$-representation $\rho_{m,t}$ is
not equivalent to $\chi\otimes\rho_{m,s}$ or $\chi\otimes\rho_{m,s}^*$ for
any character $\chi$.
\end{enumerate}
\end{prop}
\begin{remark}
Parts (a) and (b) of the proposition hold for all primitive roots of
unity $q$, not just for the specific values in the statement since the
groups in question are finite.
\end{remark}

Thanks to Theorem \ref{corresp}, up to tensor product with
a $1$-dimensional representation of $\B_m$, we can identify each
$\oW_{m,\lambda}$ with a representation of the form
$\s^2\overline{V}_{m,s}$, $\wedge^2\overline{V}_{m,s}$, or
$\overline{V}_{m,s}\otimes\overline{V}_{m,t}$.
Note that tensoring with a $1$-dimensional representation does not affect
$\bar H_{m,\lambda}$.
The image of a group in $\PGL(V)$ (resp. when $\dim V\ge 3$) is the same as its image in
$\PGL(\s^2 V)$ (resp. $\PGL(\wedge^2 V)$), since the natural homomorphisms
$\PGL(d)\to \PGL\bigl(\binom{d+1}{2}\bigr)$ and $\PGL(d)\to \PGL\bigl(\binom{d}{2}\bigr)$ are
injective for $d\ge 1$ (resp. $d\ge 3$).
To identify the closed images of $\B_m$ under the projectivized tensor products of Jones representations, we
combine Proposition \ref{tlimages} with Goursat's lemma:

\begin{lemma}[\cite{Gou}] Suppose $H\subset G_1\times G_2$ such that the compositions
$H\hookrightarrow G_1\times G_2\rightarrow G_1$ and
$H\hookrightarrow G_1\times G_2\rightarrow G_2$ are
surjective homomorphisms.  There there exist normal subgroups
$N_i\vartriangleleft G_i$ and an isomorphism $\psi:
G_1/N_1\rightarrow G_2/N_2$ such that $H$ is the graph of $\psi$,
\emph{i.e.} $(g_1,g_2)\in H$ if and only if $\psi(g_1N_1)=g_2N_2$.
\end{lemma}

Suppose $V_1$ and $V_2$ are
representation spaces of irreducible unitary representations
of $\B_m$.  Tensor product defines a natural injective map
$\PGL(V_1)\times\PGL(V_2)\to\PGL(V_1\otimes V_2)$.
If $H$, $G_1$ and $G_2$ denote the closure of the image of $\B_m$
in $\PGL(V_1\otimes V_2)$, $\PGL(V_1)$ and $\PGL(V_2)$ respectively,
then $H\hookrightarrow G_1\times G_2$ satisfies the hypotheses
of the lemma.

\begin{theorem}\label{images}
Let $\ell\ge 6$, $m\ge 1$, and $\lambda\in\Gamma(\ell)$ be such that $m-|\lambda|$
is a non-negative even integer.  Let
$$\Sigma_6 = \{[4],[4,1,1],[1,1,1,1],[2,2],[1,1],[0]\}.$$
Then,
\begin{equation}
\label{thirteen}
\bar H_{m,\lambda} \cong
\begin{cases}
\{1\}&\text{$\lambda = [m]$,} \\
\{1\}&\text{$m=2$,} \\
\{1\}&\text{$m=3$ and $\lambda = [1,1,1]$,} \\
\{1\}&\text{$m=4$ and $\lambda = [1,1,1,1]$,} \\
\PSp_{m-1}(3)&\text{$\ell=6$ and $m$ is odd,} \\
\PSp_{m-2}(3)&\text{$\ell=6$ and $\lambda\notin\Sigma_6$,} \\
\PSp_{m-2}(3)\ltimes (\Z_3)^{m-2}&\text{$\ell=6$,} \\
A_5&\text{$\ell=10$, $m=3$, and $\lambda\in\{[2,1],[1]\}$,} \\
A_5&\text{$\ell=10$, $m=4$, and $\lambda\in\{[2,2],[0]\}$,} \\
A_5\times \PSU(3)&\text{$\ell=10$, $m=4$, and $\lambda=[1,1]$,} \\
\PSU(d_{m,\frac{m-\lambda_1}{2}}) & \text{$\lambda = [\lambda_1]$,} \\
\PSU(d_{m,\frac{m-\lambda_1}{2}}) & \text{$\lambda = [\lambda_1,1,1]$,} \\
\PSU(d_{m,\frac{m}{2}}) & \text{$\lambda = [1,1,1,1]$,} \\
\end{cases}
\end{equation}
Here we employ the convention that each condition is assumed to exclude all previous ones, so that for
example the sixth case implicitly requires that $m$ is even.
In the generic case, when none of these conditions applies, we have
\begin{equation}
\label{generic}
\bar H_{m,\lambda}
=\PSU(d_{m,\frac{m-\lambda_1+\lambda_2}{2}})\times \PSU(d_{m,\frac{m-\lambda_1-\lambda_2}{2}}).
\end{equation}

\end{theorem}

\begin{proof}
To begin with, we note that the first four cases of (\ref{thirteen}) are precisely those for which
$\dim \oW_{m,\lambda} = 1$, so we may now assume
$\dim \oW_{m,\lambda} > 1$.
If $\oW_{m,\lambda}$ is the symmetric square of
some $\oV_{m,p}$, then $\bar H_{m,\lambda} \cong \bar G_{m,s}$.
Since $\binom d2>1$, the exterior square map is
injective, so the previous remark applies also in this case.  These two remarks account for
the last three cases of (\ref{thirteen}) as well as various subcases of the fifth, sixth, seventh, eighth, and ninth cases.

The remaining difficulty is to determine $\bar H_{m,\lambda}$ in the tensor product case, when we know
it is a subgroup of $\bar G_{m,s}\times \bar G_{m,t}$ mapping onto each factor.  When the two
factors are simple and non-isomorphic, then $\bar H_{m,\lambda}$ must be the whole product.
This is the situation in the tenth case of (\ref{thirteen}) and in (\ref{generic}), when
$$d_{m,\frac{m-\lambda_1+\lambda_2}{2}} \neq d_{m,\frac{m+\lambda_1+\lambda_2}{2}}.$$
If both factors are simple and isomorphic and $\bar H_{m,\lambda}$ is not the whole product,
then it must be the graph of an isomorphism.  Every automorphism of $\PSU(d)$ is either inner or the product of an inner automorphism with transpose inverse.  Therefore, if
$\bar G_{m,s}\cong\bar G_{m,t}\cong\PSU(d)$ and $\bar H_{m,\lambda}$ is the graph of an isomorphism, the representation $\rho_{m,s}$ must be equivalent (up to tensoring by a $1$-dimensional representation) to $\rho_{m,t}$ or its dual.  This is impossible by Proposition~\ref{tlimages}, so this
finishes the case that $\bar G_{m,s}$ and $\bar G_{m,t}$ are both infinite.

The only remaining cases are those where $\bar G_{m,s}$ and $\bar G_{m,t}$ are
both finite and non-trivial.  This cannot happen if $\ell=10$ (since for each $m$ there is at most one
non-trivial value of $s$ which give non-trivial finite image).  It can happen only if $\ell=6$.
Here we know \cite{BWj} that the whole image of $\B_m$ in $\overline{T}_m$ is a central extension of either $\PSp_{m-1}(3)$ or $\PSp_{m-2}(3)\ltimes (\Z_3)^{m-2}$ depending on whether $m$ is odd or even.
The tensor product of any $\oV_{m,s}$ and $\oV_{m,t}$ is contained in the symmetric square of
$\overline{T}_m$, so the image of $\B_m$ in the projectivization of any such tensor product is
a quotient of  $\PSp_{m-1}(3)$ or $\PSp_{m-2}(3)\ltimes (\Z_3)^{m-2}$ respectively.
When $m$ is odd, we therefore automatically have the fifth case of (\ref{thirteen}).
By Proposition~\ref{tlimages}, when
$m$ is even, $s=m/2-1$ gives
$\bar G_{m,s} = \PSp_{m-2}(3)\ltimes (\Z_3)^{m-2}$ and the other two values, $s=m/2$ and
$s=m/2-2$, give $\bar G_{m,s} = \PSp_{m-2}(3)$.  This means that if $s=m/2-2$, $t=m/2$,
$\bar H_{m,\lambda} = \PSp_{m-2}(3)$ (the sixth case of (\ref{thirteen})).
The remaining possibilities for $s$ and $t$ give quotients of $\PSp_{m-2}(3)\ltimes (\Z_3)^{m-2}$
which also map onto the same group, and therefore give examples belonging to the
seventh case of (\ref{thirteen}).

\end{proof}

We conclude by remarking on a striking aspect of these final cases: the tensor products of
certain pairs of irreducible representations of  $\PSp_{m-2}(3)\ltimes (\Z_3)^{m-2}$ or
$\PSp_{m-2}(3)$ turn out to be irreducible.  In particular, the two Weil representations of
$\PSp_{m-2}(3)$ have an irreducible tensor product.  It would be interesting to find other examples of
faithful projective representations which have an irreducible tensor product.  We are aware
of a number of ``sporadic'' examples but only two other infinite families, one arising from square Young diagrams in the representation theorem of $A_{n^2}$ and one from Weil representations of unitary groups over the field with two elements.

\end{document}